\documentclass[a4paper,10pt]{article}
\usepackage{amssymb,amsmath,amsthm,latexsym}
\usepackage{color}
\usepackage{textcomp}
\usepackage[active]{srcltx}
\usepackage{colortbl}

\usepackage{hyperref}

\usepackage{booktabs}

\usepackage{mathrsfs}
\usepackage{graphicx}
\usepackage{lineno}
\usepackage{adjustbox}
\usepackage{rotating}

\usepackage{longtable}

\usepackage{tikz}
\usetikzlibrary{shapes.geometric, arrows}

\tikzstyle{process} = [rectangle, minimum width=2.8cm, minimum height=0.8cm, text centered, draw=black, font=\small]
\tikzstyle{decision} = [ellipse, minimum width=2.5cm, minimum height=0.8cm, text centered, draw=black, font=\small]
\tikzstyle{arrow} = [thick,->,>=stealth]

\usepackage[multiple]{footmisc}

\usepackage{fancyhdr}

\newcommand{\keywords}[1]{%
  \par\noindent\small\textbf{\textit{Keywords---}} #1
}
\newcommand{\mscclass}[1]{%
  \par\noindent\small\textbf{MSC:} #1
}

\theoremstyle{definition}

\theoremstyle{remark}

\numberwithin{equation}{section}

\title{Algorithmic aspects of Newman polynomials and their divisors}

\author{M. Idris, J.-M. Sac-\'Ep\'ee\footnote{musbahu.idris@univ-lorraine.fr, jean-marc.sac-epee@univ-lorraine.fr, IECL, Universit\'e de Lorraine, France}}

\begin{document}

\maketitle

\newcommand{\D}{\mathbb{D}}
\newcommand{\C}{\mathbb{C}}
\newcommand{\R}{\mathbb{R}}
\newcommand{\Z}{\mathbb{Z}}
\newcommand{\dist}{\operatorname{dist}}

\renewcommand{\qedsymbol}{$\blacksquare$}

\makeatletter
\def\blfootnote{\xdef\@thefnmark{}\@footnotetext}
\makeatother

\begin{abstract}
We study the problem of determining which integer polynomials divide Newman
polynomials. First, we investigate the $8438$ known polynomials with Mahler
measure less than $1.3$. Apart from the polynomials having a positive real
root, all but three of them are shown to divide a Newman polynomial of degree
at most $1000$. Second, using an algorithmic certification method of Hare and Mossinghoff, we
exhibit several polynomials that divide no Newman polynomial. Among them is a
degree-$10$ polynomial with Mahler measure approximately equal to
$1.419404632$, which improves the best known upper bound for any universal
constant $\sigma$, if such a constant exists, governing the existence of
Newman multiples in terms of Mahler measure. Finally, writing $l(x)$ for the polynomial obtained from Lehmer's polynomial
by the substitution $x\mapsto -x$, we explicitly construct Newman polynomials
divisible by $l(x)^2$ with degrees up to $150$, and show that no Newman
polynomial of degree at most $160$ is divisible by $l(x)^3$.
\end{abstract}

\smallskip
\mscclass{11R06, 11C08, 12D10}

\section{Introduction}

Recall that the Mahler measure of a polynomial $p$ with complex coefficients is defined as
\[
M(p) := |b_0| \prod_{k=1}^{d} \max(1, |\alpha_k|),
\]
where
\[
p(x) = b_0 x^d + b_1 x^{d-1} + \dots + b_d = b_0 \prod_{k=1}^{d} (x - \alpha_k), \text{with }  b_0 \neq 0 .
\]
It is well known (\cite{BombieriVaaler1987}) that every polynomial $p(x) \in \mathbb{Z}[x]$ having Mahler measure less than $2$ must divide a polynomial with all coefficients in $\{-1, 0, 1\}$. Given this result, it is tempting to ask the following question: does this result have any chance of persisting if we impose stricter conditions on the coefficients, namely that they are no longer in $\{-1, 0, 1\}$, but only in $\{0, 1\}$?

Recall that a \emph{Newman polynomial} is a polynomial whose coefficients are either $0$ or $1$, with constant term $1$. In the remainder of this article, we denote by $\mathcal{N}$ the set of all Newman polynomials. In fact, we can immediately refine the question above through a few simple observations.

First, a polynomial with a positive real root cannot divide a Newman polynomial. More generally, one should exclude nonnegative real roots, since a zero root is also impossible for a divisor of a Newman polynomial, whose constant term is equal to $1$. In the computations below, however, the polynomials under consideration have nonzero constant term; therefore this obstruction will usually be stated in the equivalent form of having a positive real root.

A second factor stems from the result of A. M. Odlyzko and B. Poonen (\cite{OdlyzkoPoonen1993}), which states that any nonzero root $\alpha$ of a Newman polynomial satisfies $1/g<\vert \alpha\vert< g$, where $g \approx 1.6180 $ is the golden ratio. This result might lead us to think that an integer-coefficient polynomial without a nonnegative real root, provided that its Mahler measure is smaller than $g$, could divide a Newman polynomial. In this regard, a paper by K. G. Hare and M. J. Mossinghoff (\cite{MossinghoffHare2014}) proves that the golden ratio is a natural bound for a wide range  of classes of polynomials. Recall that a real algebraic integer strictly greater than $1$ is a \emph{Salem number} if all its conjugates (that is, the other roots of its minimal polynomial) have modulus less than or equal to $1$, and at least one conjugate has modulus exactly $1$. Recall also that a \emph{Pisot number} is a real algebraic integer strictly greater than $1$ whose conjugates all have modulus strictly less than $1$. Given these definitions, we say that a real number $\alpha$ is a \emph{negative Pisot number} (resp. a \emph{negative Salem number})  if $-\alpha$ is a Pisot number (resp. a Salem number). In \cite{MossinghoffHare2014}, the authors prove that every negative Pisot number within the interval $(-g,-1)$ that has no positive conjugates, as well as every negative Salem number in $(-g,-1)$ obtained through Salem's construction applied to small negative Pisot numbers, satisfies a Newman polynomial. Additionally, they identify all negative Salem numbers within $(-g,-1)$ with a degree of at most $20$ and confirm that each of them is a root of a Newman polynomial.

Another interesting result is that of A. Dubickas (\cite{Dubickas2003}), which states that any product of cyclotomic polynomials that does not contain powers of $x-1$ necessarily divides a Newman polynomial. Related questions for sparse polynomials with restricted coefficients have also been studied by D. Dutykh and J.-L. Verger-Gaugry. In \cite{DutykhVerger-Gaugry2018}, they investigate almost Newman lacunary polynomials, with particular emphasis on reducibility, cyclotomic and noncyclotomic factors, reciprocal and nonreciprocal parts, and the geometry of their zero sets. Further connections between restricted coefficient polynomials, algebraic bases, rewriting trails, ultimately periodic representations, almost Newman polynomials, and Lehmer-type questions are developed in \cite{DutykhVerger-Gaugry2021}.

However, it turns out that the golden ratio is a bound that is too large, according to the same article by K. G. Hare and M. J. Mossinghoff (\cite{MossinghoffHare2014}), where the authors exhibit five polynomials whose Mahler measure is smaller than $g$, yet none of them divides a Newman polynomial. The smallest of the five Mahler measures is approximately equal to $1.55601$, which corresponds to the polynomial \(x^6 - x^5 - x^3 + x^2 + 1\).

In light of all this, the natural question that arises is the following: \emph{Does there exist a constant $\sigma$ such that if the integer-coefficient polynomial $p$ has no nonnegative real root and has a Mahler measure smaller than $\sigma$, then $p$ necessarily divides a Newman polynomial?}

Naturally, the result of K. G. Hare and M. J. Mossinghoff (\cite{MossinghoffHare2014}) implies that $\sigma$, if it exists, is at most approximately $1.55601$. This upper bound was improved by P. Drungilas, J. Jankauskas, and J. \v{S}iurys in \cite{DrungilasJankauskasSiurys}, who studied the existence of Littlewood and Newman multiples of Borwein polynomials by means of a decision algorithm for monic integer polynomials with no roots on the unit circle. As part of this work, they produced new examples of polynomials with no Newman multiple and showed that $\sigma$, if it exists, must in fact be less than $1.436632261$.

We first chose to investigate polynomials with small Mahler measure, since such polynomials have been studied extensively for decades. In particular, tables of such polynomials are available, some of which are certified as complete. More precisely, we will examine the list of all known polynomials having a Mahler measure less than $1.3$ and a degree less than or equal to $180$, with the aim of determining whether the polynomials in this list divide any Newman polynomials. This list was available for many years on a website maintained by M. J. Mossinghoff. Since that website is no longer available, we use the copy reproduced in the file \texttt{Known180.txt} included in our dataset \cite{ComputationalResults}. For clarity and readability, we refer to this list throughout the article as \emph{Known180}. This non-exhaustive list of $8438$  polynomials is complete up to degree $44$ (\cite{MossinghoffRhinWu2008}). In the following section, we will present an algorithm designed to explicitly provide a Newman polynomial (if it exists) that a polynomial from this list divides, provided that the sought Newman polynomial does not have a degree strictly greater than 1000. Apart from polynomials that have a positive real root, we will see that only three polynomials remain unresolved, in the sense that no Newman multiple of degree at most \(1000\) was found for them.

\section{The case of polynomials with small Mahler measure}
As announced, in this first approach we will focus on polynomials in $\mathbb{Z}[x]$ whose Mahler measure is strictly less than $1.3$. All polynomials in \emph{Known180} are monic and irreducible. Moreover, they are reciprocal by a classical result of C. J. Smyth proved in \cite{Smyth1971}. We have developed and implemented an algorithm designed to take as input each polynomial from \emph{Known180} and search (up to degree $1000$), for a Newman polynomial that has this input polynomial as one of its factors. In the next subsection, we will precisely explain how this algorithm works.
\subsection{Description of the algorithm} \label{Description of the algorithm}
An input polynomial, chosen from \emph{Known180},  is a reciprocal polynomial of even degree (say $2d$). It can be written as
\[
p(x) = x^{2d} + a_1 x^{2d-1} + a_2 x^{2d-2} + \ldots + a_{d-1} x^{d+1} + a_d x^d + a_{d-1} x^{d-1} + \ldots + a_2 x^{2} + a_1 x + 1,
\]
where the $a_i$'s are explicitly given. If $p$ is not already a Newman polynomial, and if it does not have any strictly positive real root (which would prevent it from dividing a Newman polynomial), we search for a polynomial $q(x) = x^m + b_{m-1}x^{m-1} + \ldots + b_1 x + 1 \in \mathbb{Z}[x]$ such that $p(x)q(x)$ is a Newman polynomial. The constant term of $q$ is forced to be $1$, since the input polynomial $p$ has constant term $1$ and the product $p(x)q(x)$ must also have constant term $1$. Expanding the product \( p(x)q(x) \) and identifying the coefficients $c_k (k = 0, 1, \dots, 2d+m)$ of each power of \( x \) in this product (given that the $a_i$'s are explicitly known, it follows that each $c_k$ is a linear expression in $b_1, \dots, b_{m-1}$), this comes down to solving the system of linear inequalities with integer unknowns $b_1, \dots, b_{m-1}$:
\begin{equation} \label{eq:systemConstraints}
0 \leq c_k \leq 1, \quad \text{for } k = 0, 1, \dots, 2d+m.
\end{equation}
We increase the value of $m$ starting from $1$ until we find a solution (that is, a polynomial $q$) satisfying the system of constraints \eqref{eq:systemConstraints}, in which case $p(x)$ is indeed a factor of the Newman polynomial $p(x)q(x)$. We do not allow $m$ to keep increasing if $2d + m$ becomes greater than $1000$. As an example, let us suppose we want to determine whether the polynomial $p(x) = x^{16} + x^{15} - x^{11} - x^8 - x^5 + x + 1$ from \emph{Known180} is a factor of a Newman polynomial. For $m$ from $1$ to $21$, system \eqref{eq:systemConstraints} has no solution. For $m = 22$, we search for integers $b_1, \ldots, b_{21}$, that is, a polynomial $q(x) = x^{22} + b_{21}x^{21} + \ldots + b_1x + 1$ (not necessarily reciprocal), such that $p(x)q(x)$ is a Newman polynomial. In this particular case, System \eqref{eq:systemConstraints} becomes

Find $b_1$, $b_2$, \ldots, $b_{21}$ in $\mathbb{Z}$ such that
\begin{equation} \label{eq:exampleSystem}
\left\{
\begin{array}{l}
0 \leq b_{21} + 1 \leq 1 \\
0 \leq b_{21} + b_{20} \leq 1 \\
0 \leq b_{20} + b_{19} \leq 1 \\
0 \leq b_{19} + b_{18} \leq 1 \\
\vdots \\
0 \leq b_{5} + b_{4} - 1 \leq 1 \\
0 \leq b_{4} + b_{3} \leq 1 \\
0 \leq b_{3} + b_{2} \leq 1 \\
0 \leq b_{2} + b_{1} \leq 1 \\
0 \leq b_{1} + 1 \leq 1
\end{array}
\right.
\end{equation}
Solving the system with integer unknowns \eqref{eq:exampleSystem} yields the polynomial $q(x) = x^{22} + x^{20} + x^{18} + x^{15} + x^{13} + 2x^{11} + x^9 + x^7 + x^4 + x^2 + 1$, and it is easy to verify that \( p(x)q(x) = x^{38} + x^{37} + x^{36} + x^{35} + x^{34} + x^{27} + x^{11} + x^4 + x^3 + x^2 + x + 1\), which is indeed a Newman polynomial. We solve the system of inequalities \eqref{eq:exampleSystem} using the mixed-integer linear programming package Gurobi, which is called from a program written in Julia. We will provide further technical details later in the article. We did this for the $8438$ polynomials in \emph{Known180}.

\subsection{Calculation results}
We have recorded the results of the computations in the accompanying dataset \cite{ComputationalResults}. Unlike the result files described below, the file \texttt{Known180.txt}, which contains the list \emph{Known180}, records only the nonredundant coefficients of each polynomial, since all polynomials in \emph{Known180} are reciprocal. The other files contain all the coefficients of each polynomial.

The file \texttt{newman.txt} contains all the polynomials from \emph{Known180} that are already Newman polynomials. Each line contains the polynomial along with its degree, Mahler measure, and the full list of coefficients.

The file \texttt{positive.txt} contains all the polynomials from \emph{Known180} that, having a positive real root, cannot divide a Newman polynomial. Each line contains the polynomial with its degree, Mahler measure, one of its positive real roots, and the list of its coefficients.

The last three files, namely \texttt{outputP.txt}, \texttt{outputQ.txt}, and \texttt{outputD.txt}, are interconnected. On each line $k$ of the file \texttt{outputP.txt}, there is a polynomial $p$ from \emph{Known180}, with its degree, Mahler measure, and list of coefficients. On the same line $k$ of \texttt{outputQ.txt}, there is a polynomial $q$ such that the product $pq$ is a Newman polynomial. Line $k$ of \texttt{outputQ.txt} contains the degree of $q$ and its list of coefficients. Finally, each line $k$ of the file \texttt{outputD.txt} contains the product of the polynomials $p$ and $q$ found on lines $k$ of \texttt{outputP.txt} and \texttt{outputQ.txt}. The file \texttt{outputD.txt} therefore contains only Newman polynomials.

To verify that the Newman multiples are correct, the reader may use the short script \texttt{verifNewman.gp} (available at \cite{ComputationalResults}) written in GP/PARI, which should be placed in the directory where the files \texttt{outputP.txt}, \texttt{outputQ.txt}, and \texttt{outputD.txt} have previously been downloaded. One simply needs to indicate, in the first line of the script, which line $k$ of the files \texttt{outputP.txt}, \texttt{outputQ.txt}, and \texttt{outputD.txt} is to be checked. The script displays the polynomial $p$ appearing on line $k$ of \texttt{outputP.txt} (the degree and the Mahler measure are skipped), the polynomial $q$ appearing on line $k$ of \texttt{outputQ.txt} (the degree is skipped), computes and displays the product $pq$, and then displays the Newman polynomial $d$ appearing on line $k$ of \texttt{outputD.txt} (the degree is skipped). The script confirms that one indeed has $pq = d$.

It is worth noting that among the $8438$ polynomials examined, if one excludes those that were discarded from the outset because they have a positive real root, only three polynomials remain unresolved, and it is unknown whether they divide a Newman polynomial. 

The only three polynomials in \emph{Known180} for which we were unable to exhibit multiples in the set $\mathcal{N}$ of Newman polynomials are recorded in the file \texttt{nonsol.txt} in the dataset \cite{ComputationalResults}. Our algorithm certifies only that these polynomials do not divide any Newman polynomial of degree at most \(1000\). It remains open whether they may occur as factors of Newman polynomials of higher degree.

In the following section, we examine a method which, under certain conditions, can certify that a given polynomial divides no Newman polynomial.

\section{Polynomials that divide no Newman polynomial}
In this section, we employ an algorithmic method introduced in \cite{MossinghoffHare2014} by K. G. Hare and M. J. Mossinghoff to construct examples of polynomials with integer coefficients that do not divide any Newman polynomial. We begin by briefly outlining the underlying principles of this algorithm. For the sake of clarity, we strictly adhere to the notations, particularly the names of the sets introduced in \cite{MossinghoffHare2014}.

Let $\beta$ be an algebraic integer with $|\beta| > 1$, and define the closed disk
\[
I'(\beta) = \left\{ z \in \mathbb{C} : |z| \leq \frac{|\beta|}{|\beta| - 1} \right\}.
\]

We iteratively build the set of values
\[
\mathcal{N}'(\beta, d) = \{F(\beta) : F \in \mathcal{N}, \deg(F) \leq d \} \cap I'(\beta),
\]
and define the \textit{total reachable} set as
\[
\mathcal{N}'(\beta) = \bigcup_{d \geq 0} \mathcal{N}'(\beta, d).
\]
The set $\mathcal{N}'(\beta, d)$ can be computed recursively: starting with $\mathcal{N}'(\beta, 0) = \{1\}$, each subsequent set is generated by the recurrence
\[
\mathcal{N}'(\beta, d+1) = \mathcal{N}'(\beta, d) \cup \left( \{ \beta \omega : \omega \in \mathcal{N}'(\beta, d) \} \cap I'(\beta) \right) \cup \left( \{ \beta \omega + 1 : \omega \in \mathcal{N}'(\beta, d) \} \cap I'(\beta) \right).
\]
This construction exploits the fact that for any $F \in \mathcal{N}$, if $F(\beta) \notin I'(\beta)$, then neither $\beta F(\beta)$ nor $\beta F(\beta) + 1$ can lie in $I'(\beta)$.

Given that we are seeking a polynomial that divides no Newman polynomial, there is no point in continuing the construction of \( \mathcal{N}'(\beta, d) \) if $0 \in \mathcal{N}'(\beta, d)$  for some $d$, since this implies that $\beta$ satisfies a certain Newman polynomial \( F \), and therefore the minimal polynomial of $\beta$ divides this Newman polynomial $F$. The case of interest to us arises when, for some \( d \), it turns out that $\mathcal{N}'(\beta, d+1) = \mathcal{N}'(\beta, d)$, and simultaneously $0 \notin \mathcal{N}'(\beta, d)$. In this case, it becomes clear that $\mathcal{N}'(\beta)$ is finite and that $0 \notin \mathcal{N}'(\beta)$, so that \( \beta \) is a root of no Newman polynomial, and therefore the minimal polynomial of \( \beta \) divides no Newman polynomial.

In \cite{MossinghoffHare2014}, K.~G.~Hare and M.~J.~Mossinghoff provide a list of monic polynomials with integer coefficients, all of whose roots lie in the slit annulus
\[
\{z \in \mathbb{C} : g^{-1} < |z| < g\} \backslash \mathbb{R}^+,
\]
and prove, using the algorithmic method described above, that these polynomials do not divide any Newman polynomial. Among these polynomials, the one with the smallest Mahler measure is the polynomial $x^6 - x^5 - x^3 + x^2 + 1$, whose Mahler measure is approximately $1.556014485$. The significance of their result lies in establishing that any constant
$\sigma$ such that every integer-coefficient polynomial of Mahler measure
less than $\sigma$ and with no nonnegative real root divides a Newman
polynomial, if such a constant exists, must satisfy
\[
\sigma \leq M(x^6-x^5-x^3+x^2+1)\approx 1.556014485.
\] 
This bound was improved by P. Drungilas, J. Jankauskas and J. \v{S}iurys in \cite{DrungilasJankauskasSiurys}, where the authors enlarged the list given in \cite{MossinghoffHare2014} by adding several polynomials of Mahler measure less than \(1.556014485\), one of which has Mahler measure approximately equal to \(1.436632261\). Consequently, any such constant \(\sigma\), if it exists, must be no greater than this Mahler measure.

Our aim is now to further improve the bound established in
\cite{DrungilasJankauskasSiurys}. Along the way, we obtain a list of
certified examples of degrees $8$, $9$, $10$, $11$, and $12$, with Mahler
measure not exceeding $1.556014485$; some of these recover previously known
examples, while others appear to be new. In particular, one degree-$10$
polynomial has Mahler measure approximately equal to $1.419404632$,
thereby improving the bound of \cite{DrungilasJankauskasSiurys}. We proceed in two stages. First, we identify candidate polynomials, which we hope do not divide any Newman polynomial, by rapidly testing, with the aid of the algorithm described in Section~2, all monic polynomials of degrees $6$, $7$, $8$, $9$, $10$, $11$, and $12$ whose coefficients lie between \(-2\) and \(2\). In view of the high speed of this preliminary screening step, we did not attempt to reduce the amount of computation by exploiting symmetries or other similar reductions. This step took fewer than three hours on a laptop running Ubuntu 24.04.4 LTS, equipped with an Intel(R) Core(TM) Ultra 9 185H processor (16 cores: 6 performance, 8 efficient, and 2 low-power efficient cores; 22 threads) and 64 GB of installed RAM (about 62 GB usable under Linux). In brief, for each tested polynomial, we solve systems of inequalities with integer unknowns in order to test whether the candidate polynomial divides some Newman polynomial of degree at most $200$. The algorithm used in this first stage has been fully implemented in the Julia programming language (\textsc{Julia}~1.11.6). The optimization problem is formulated within the framework of the JuMP modeling language, while the search for a solution that respects the constraints imposed is carried out using the Gurobi library (Gurobi Optimizer~12.0.3) via a wrapper maintained by the JuMP community. Gurobi was run with its default parameters. Each candidate polynomial found at the preselection stage is then verified using the algorithm from \cite{MossinghoffHare2014}, and is either retained or discarded according to the outcome of this verification. The preselection step produced 38 irreducible polynomials, all of whose roots lie in the slit annulus $\{z \in \mathbb{C} : g^{-1} < \vert z \vert < g\} \backslash \mathbb{R}^+$, and none of which divides any Newman polynomial up to degree 200. When several polynomials had the same Mahler measure, only one representative was retained. The verification step certified 33 of them; the remaining 5 were left undecided by the computation, since the number of elements in the corresponding set does not stabilize within the maximum allowed depth (set to 200), while the maximum allowed cardinality of the set was fixed at 40,000,000 in order to avoid excessive memory consumption.

After performing all computations and verifications, the resulting polynomials are listed in Table~\ref{tab:mahler-polynomials-full}, which also specifies, for each polynomial, the chosen root $\beta$ used in the computation, the stabilization depth that was obtained, and the cardinality of the corresponding set.

{\scriptsize
\begin{longtable}{p{0.34\textwidth}p{0.13\textwidth}p{0.18\textwidth}p{0.10\textwidth}p{0.10\textwidth}}
\caption{Certified examples of polynomials with Mahler measure not exceeding $1.556014485$ and not dividing a Newman polynomial}
\label{tab:mahler-polynomials-full}\\
\toprule
\textbf{Polynomial} & \textbf{Mahler measure} & \textbf{Chosen conjugate $\beta$} & \textbf{$d_{\mathrm{stab}}$} & \textbf{$|\mathcal{N}'(\beta)|$} \\
\midrule
\endfirsthead

\endhead

\endfoot

\bottomrule
\endlastfoot

$x^{10} - x^8 - x^5 + x + 1$ & 1.419404632 & $1.1307-0.1621i$ & 17 & 768 \\
$x^9 - x^6 - x^5 - x^4 + x^2 + x + 1$ & 1.436632261 & $1.1403-0.1518i$ & 14 & 433 \\
$x^{12} - x^{11} - x^6 + x^5 + 1$ & 1.448290492 & $1.1426-0.1884i$ & 15 & 388 \\
$x^{10} - x^7 - x^6 - x^5 + x^2 + x + 1$ & 1.475517312 & $1.1562-0.0996i$ & 13 & 232 \\
$x^{10} - x^9 - x^5 + x^4 + 1$ & 1.475928627 & $1.1569-0.2252i$ & 12 & 190 \\
$x^{11} - x^9 - x^8 + x + 1$ & 1.477652735 & $1.1928-0.0507i$ & 10 & 82 \\
$x^{10} - x^8 - x^7 + x^5 - x^3 + x + 1$ & 1.481909009 & $1.1576-0.1629i$ & 12 & 200 \\
$x^9 - x^7 - x^6 + x + 1$ & 1.483444878 & $1.1860-0.1451i$ & 10 & 85 \\
$x^{12} - x^{11} - x^8 + x^7 + 1$ & 1.488053700 & $1.1620-0.1901i$ & 11 & 165 \\
$x^8 - x^7 - x^4 + x^3 + 1$ & 1.489581321 & $1.1736-0.2801i$ & 10 & 95 \\
$x^{12} - x^9 - x^7 - x^5 + x^4 + x^2 + 1$ & 1.503425675 & $1.1164-0.1208i$ & 19 & 1932 \\
$x^{12} + x^{11} - x^9 - x^8 - x^7 - x^6 - x^5 + x^3 + x^2 + x + 1$ & 1.504539808 & $1.0939-0.0894i$ & 27 & 15322 \\
$x^{12} - x^{10} - x^7 - x^6 + x^4 + x + 1$ & 1.505811890 & $1.1615-0.0544i$ & 13 & 196 \\
$x^{12} - x^9 - x^7 - x^5 + x^2 + x + 1$ & 1.515809085 & $1.1305-0.0819i$ & 17 & 753 \\
$x^8 - x^6 - x^5 + x + 1$ & 1.518690904 & $1.1765-0.1929i$ & 10 & 102 \\
$x^{10} - x^9 - x^7 + x^6 + x^4 - x^3 + 1$ & 1.519065984 & $1.1666-0.2719i$ & 11 & 120 \\
$x^{12} - x^{10} - x^7 + x + 1$ & 1.519151476 & $1.1402-0.1065i$ & 15 & 451 \\
$x^{12} - x^8 - x^7 - x^6 - x^5 + x^3 + x^2 + x + 1$ & 1.520043072 & $1.1339-0.0516i$ & 16 & 637 \\
$x^{12} - x^{10} - x^9 + x^7 + 1$ & 1.521929346 & $1.1381-0.1765i$ & 15 & 491 \\
$x^{12} - x^{11} - x^{10} + x^9 + x^8 - x^6 + x^4 + x^3 - x^2 - x + 1$ & 1.522096479 & $1.1440-0.4619i$ & 17 & 320 \\
$x^{12} - 2x^{11} + x^{10} + x^9 - x^8 + x^7 - x^6 + x^5 - x^4 + x^3 + x^2 - 2x + 1$ & 1.526433287 & $1.0877-0.5859i$ & 28 & 1325 \\
$x^{12} - x^{10} - x^7 + x^3 + 1$ & 1.531401303 & $1.1328-0.1357i$ & 16 & 655 \\
$x^{11} - x^9 - x^8 + x^6 + x^5 - x^4 - x^3 + x + 1$ & 1.534099819 & $1.1333-0.1744i$ & 17 & 664 \\
$x^{11} - x^8 - x^6 + x + 1$ & 1.534872750 & $1.1048-0.1448i$ & 53 & 89226 \\
$x^8 - x^6 - x^5 - x^4 + x^2 + x + 1$ & 1.536566472 & $1.2220-0.1122i$ & 8 & 41 \\
$x^9 - x^8 - x^6 + x^5 + 1$ & 1.536913983 & $1.1920-0.2548i$ & 8 & 60 \\
$x^{11} - x^{10} - x^9 + x^8 + x^7 - x^6 - x^5 + x^4 + x^3 - x^2 + 1$ & 1.537810676 & $1.1679-0.2533i$ & 11 & 123 \\
$x^{11} - x^7 - x^6 + x^3 + x^2 + x + 1$ & 1.539383026 & $1.1336-0.2531i$ & 33 & 5293 \\
$x^{11} - x^{10} + x^3 + 1$ & 1.542629108 & $1.1786-0.3131i$ & 10 & 75 \\
$x^{12} - x^{11} - x^9 + x^7 + x^5 - x^4 + 1$ & 1.547469802 & $1.2124-0.1812i$ & 8 & 45 \\
$x^9 - x^7 - x^6 + x^4 + 1$ & 1.550687063 & $1.1423-0.2286i$ & 22 & 454 \\
$x^{11} - x^{10} - x^5 + x^4 + 1$ & 1.553591690 & $1.1450-0.2051i$ & 15 & 348 \\
$x^6 - x^5 - x^3 + x^2 + 1$ & 1.556014485 & $1.1910-0.3710i$ & 8 & 44 \\
\end{longtable}
}
The certification procedure used to prove that a polynomial divides no Newman polynomial was implemented independently in two different environments: first in Julia, and then in GP/PARI (version 2.15.4). In the Julia implementation, polynomial manipulations were carried out using the package Polynomials.jl, while irreducibility tests were performed with Nemo.jl. In GP/PARI, the corresponding computations were reproduced using its native polynomial routines. The GP/PARI implementation is slower than the Julia version, but it provided an independent cross-check of the entries reported in Table~\ref{tab:mahler-polynomials-full}. We note that this table is not certified to be exhaustive, since 5 polynomials were left undecided by the certification procedure.

In the following section, we investigate Newman polynomials from a different perspective, focusing on their ability to have repeated zeros outside the unit disk. At the same time, this second application aims to highlight the flexibility of the MILP framework developed in this article by showing that it can also be adapted to the construction of Newman polynomials divisible by prescribed repeated noncyclotomic factors.

\section{Newman polynomials with repeated zeros}
\subsection{State of the art}

In \cite{BorweinErdelyi1996}, P.~Borwein and T.~Erd\'elyi recall that, after observing through numerical computations that no Newman polynomial of degree $\leq 25$ has a repeated root of modulus strictly greater than $1$, A.~Odlyzko raised the question of whether every Newman polynomial necessarily has its repeated roots either at $0$ or on the unit circle. This question was resolved in \cite{Mossinghoff2003} by M.~J.~Mossinghoff.

Following Mossinghoff's notation, we set
\[
l(x) = x^{10} - x^9 + x^7 - x^6 + x^5 - x^4 + x^3 - x + 1,
\]
which is obtained from Lehmer's polynomial by the substitution $x\mapsto -x$.
We then define $g(x)=l(x)^2$. Mossinghoff carried out an exhaustive
search up to high degrees in order to identify Newman polynomials
divisible by $l(x)^2$, the interest being that $l(x)$ has a real root of
modulus strictly greater than $1$. In the first part of his article, he
examined all Newman polynomials up to degree $60$, discovering two
suitable examples, one of degree $59$ and one of degree $60$. He then
extended his search by focusing on reciprocal Newman polynomials up to
degree $100$, leading to the discovery of $21$ additional Newman
polynomials divisible by $l(x)^2$. Six further examples were obtained
using alternative choices for $g(x)$, the largest degree among these
being $105$.

Mossinghoff's method, which is highly efficient, consists in encoding each Newman polynomial $f$ as an integer $w = f(2)$. Prior to performing any polynomial division, a large number of candidates are discarded by imposing the necessary congruence conditions
\[
f(2) \equiv 0 \pmod{g(2)}
\qquad\text{and}\qquad
f(-2) \equiv 0 \pmod{g(-2)}.
\]
This substantial arithmetic filtering reduces the search space to a small subset of cases, for which exact polynomial divisibility then needs to be tested.

An alternative approach, if one wishes to avoid exhaustive searches—which remain time-consuming even when the number of computations is reduced—is to apply our integer optimization method by treating the polynomial $l(x)^2$ in the same manner as we did for all polynomials in the list \emph{Known180}. In other words, it suffices to search for a polynomial $q(x)$ with integer coefficients such that $q(x)\,l(x)^2$ is a Newman polynomial. As in Section~\ref{Description of the algorithm}, this immediately reduces to solving a system of integer inequalities.

To carry out these computations efficiently, we opted for an accelerated alternative version of our algorithm, in which the system of inequalities is replaced by a small number of linear equations. In practice, this approach allowed us to handle high degrees directly (we rapidly reached degree $150$) while also removing the restriction to reciprocal Newman polynomials. The presentation of this method is the subject of the following section.

\subsection{Algorithm for cancelling the remainder of the Euclidean division}
\subsubsection{Outline of the algorithm}
As above, we therefore consider the polynomial
\[
l(x)
  \;=\;
  x^{10} - x^9 + x^7 - x^6 + x^5 - x^4 + x^3 - x + 1,
\]
and we set $g(x) = l(x)^2$.
For an integer $d \ge \deg(g)$, we seek a Newman polynomial of degree~$d$, that is, a polynomial of the form
\[
N(x)
 \;=\;
 x^{d} + a_{d-1} x^{d-1} + \cdots + a_{1} x + 1,
 \qquad a_1, \ldots, a_{d-1} \in \{0,1\},
\]
which is a multiple of $g(x)$.

The Euclidean division of $N(x)$ by $g(x)$ can be written as
\[
N(x) \;=\; g(x)\,H(x) + R(x),
\qquad \text{with } \deg(R) < \deg(g).
\]
Clearly,
\[
g(x)\mid N(x)
\quad \text{is equivalent to} \quad
R(x)=0.
\]

For a fixed degree $d$, the coefficients of the remainder $R(x)$ are simple linear expressions in the binary variables $a_1,\dots,a_{d-1}$, with the leading and constant
coefficients of $N(x)$ fixed equal to $1$. Consequently, the divisibility condition $g(x)\mid N(x)$ reduces to the simultaneous vanishing of these linear expressions. As a result, the search for a multiple of $l(x)^2$ among Newman polynomials of degree $d$ amounts to solving a system of linear equations consisting of $\deg(g)$ equations in the $d-1$ binary unknowns $a_1,\ldots,a_{d-1}$.

To solve this system, we employ an algorithm of the same nature as that used in Section~\ref{Description of the algorithm}, with the difference that Section~\ref{Description of the algorithm} deals with a system of linear inequalities, whereas the present setting involves a system of linear equations. In both cases, the core of the method relies on a mixed-integer linear programming (MILP) solver based on the branch-and-cut technique. These computations were carried out in the same software and hardware environment as described in Section~3, again using Gurobi with its default parameters.

\subsubsection{Computational Results and Comments}

We have chosen to present one polynomial per degree (whenever a suitable polynomial exists for the degree under consideration), up to degree $150$. As in \cite{Mossinghoff2003}, Newman polynomials are represented in Table~\ref{tablehexa} using hexadecimal encoding. Since a Newman polynomial has coefficients in $\{0, 1\}$, it may be identified with the binary word formed by its coefficients, read from the highest-degree term down to the constant term. Equivalently, this binary word is the base-2 expansion of the integer $f(2)$. The hexadecimal notation used in Table~\ref{tablehexa} is therefore nothing more than a compact representation of this binary encoding, with each hexadecimal digit encoding four consecutive binary digits. We recover several polynomials that already appear in Table~1 of \cite{Mossinghoff2003}, but new polynomials also arise well before exceeding degree~$100$, which is where Table~1 of \cite{Mossinghoff2003} ends, since our polynomials are not required to be reciprocal. Thus, for degree~$93$, we provide the non-reciprocal polynomial
\[
\begin{array}{l}
x^{93} + x^{92} + x^{86} + x^{84} + x^{83} + x^{81} + x^{80} + x^{79}
+ x^{75} + x^{73} + x^{72} + x^{71} + x^{70} + x^{69} + x^{66} +\\
x^{64} + x^{63} + x^{61} + x^{60} + x^{59} + x^{57} + x^{55}
+ x^{53} + x^{52} + x^{51} + x^{48} + x^{47} + x^{45} + x^{42} + x^{40} +\\
x^{39} + x^{37} + x^{34} + x^{32} + x^{28} + x^{27} + x^{25}
+ x^{24} + x^{22} + x^{19} + x^{16} + x^{14} + x^{13} + x^{12} + x^{11} +\\
x^{8} + x^{5} + x + 1
\end{array}
\]
which is indeed a multiple of $l(x)^2$.

{\normalfont
\begin{longtable}{r l@{\hspace{0.6em}}r l}
\caption{Newman polynomials divisible by \(l(x)^2\), encoded in hexadecimal representation}\label{tablehexa}\\
\textbf{deg($f$)} & \textbf{Hexadecimal encoding of $f$}
& \textbf{deg($f$)} & \textbf{Hexadecimal encoding of $f$} \\
\hline
\endfirsthead

\endhead

59  & \texttt{C49E23C93C47923} & 119 & \texttt{C4FAD7BA97C2CFD6B94237A41AD7FB} \\
60  & \texttt{1B33F1364D91F99B} & 120 & \texttt{18995E6FE1AA04687E78537E4B41F9B} \\
79  & \texttt{C42E67FE42427FE67423} & 121 & \texttt{37EF7A677A9AB5E1D59EEC9B8C7DE7B} \\
84  & \texttt{1BF66C43EC4446F846CDFB} & 122 & \texttt{60B7736745938332D21C93A79D37923} \\
85  & \texttt{37EDB84FDB61B6FC876DFB} & 123 & \texttt{C7CE2C0A5A59DD516336EB1318D7CE3} \\
90  & \texttt{639D99136318C63644CDCE3} & 125 & \texttt{30468A627162BA07CCECF121B0A75323} \\
92  & \texttt{1A59E126B6D4A56DAC90F34B} & 126 & \texttt{6CCFF45B93FFFE6675739863EA735323} \\
93  & \texttt{305B8BE5BAB9A5A51B497923} & 127 & \texttt{C73EB19A97352444B67C10D1D4C556A9} \\
94  & \texttt{694A44C6EF9B7A006117F99B} & 128 & \texttt{189F5AF154AAD290602D01EE4DC52F723} \\
95  & \texttt{C49E22400624F1364D91F99B} & 129 & \texttt{342AD5C57853DAA3E09EDFD732A54E94B} \\
96  & \texttt{182345F5AD226C896B5F45883} & 130 & \texttt{690766AEDA956835DF11B2276465572B3} \\
97  & \texttt{3722F1D73627FF91B3AE3D13B} & 131 & \texttt{C49ECA4E15B08C277234B3BA5C8F77F9B} \\
99  & \texttt{DFB6E49A554BFFE5EFAE3D13B} & 132 & \texttt{1A5FE9F11B73F71F0AF05B5EBC6141F99B} \\
100 & \texttt{1BC3DC3ED48DFBF6256F87787B} & 133 & \texttt{3666488BB4F195FCFF9B85E4C71A07787B} \\
101 & \texttt{34A5C8AF55736DB3AABD44E94B} & 134 & \texttt{69D7AED27383B4060B55FD28E957C5FDAF} \\
102 & \texttt{691F21E23252CA266365057ADB} & 135 & \texttt{F5A72ACE882379E14F5AD706C1904C1923} \\
103 & \texttt{DCBB435B62E0BD0746DAC2DD3B} & 136 & \texttt{1A746E02E27973B9E2480F7867E7BB235CB} \\
104 & \texttt{1B27583B6960DA5F53B9982DD3B} & 137 & \texttt{310BFBDF1989D599740B3B4917675C67FE3} \\
105 & \texttt{3787883F743027DA9D6EACA5DCB} & 138 & \texttt{694B24458A39ABB9D401777F80CF5E6694B} \\
106 & \texttt{6217333EE5D5EABD5D3BE667423} & 139 & \texttt{F77EABAFE546557B39FD2FE643A57D572B3} \\
107 & \texttt{C49E225B33E55E6B60AC954F99B} & 140 & \texttt{1BDBF94508CC7F61534EE99497E25C52F723} \\
108 & \texttt{1A52CC31D2106C65C0A298414923} & 141 & \texttt{319396AAA2ACEFA97B570D49D3711992FF4B} \\
109 & \texttt{3667E26C9B23FF7FE23C93C47923} & 142 & \texttt{6905A4AE7133EA8ECAD9EB9EF97FA7C47923} \\
110 & \texttt{62495296F6CD2243B8C1DB5F5883} & 143 & \texttt{D23F21F4DD494583F0CE0456FC0701F97923} \\
111 & \texttt{DB5EA0A6DEF59751FF6DD407FB7B} & 144 & \texttt{1885C9D6A84A3F250955F671BDACEAD727423} \\
112 & \texttt{1885CC23413B3AB8CD03B39E95323} & 145 & \texttt{34BF60CF7EEF23F169BC3C0DEA8515D0FA5FB} \\
113 & \texttt{3667EF3AF08B1450C892ACDDC7923} & 146 & \texttt{6E5CC38A712EDCA345F90B3F1ED96FDD37683} \\
114 & \texttt{621730B5A7FFBD71616397953FB7B} & 147 & \texttt{C49E87DABEB95ED9DC8182996C4E23ED5F923} \\
115 & \texttt{C10F28AC0AC78958180C6DC97BDFB} & 148 & \texttt{189684517CAAD4D9863DFA9CA16CDCFC03F1FB} \\
116 & \texttt{1BF1C8850D11BFDA96DEB370352123} & 149 & \texttt{36428A55D5405C555FA57D0AD1ED6B2011294B} \\
117 & \texttt{31212827E3D4D487E93F27FFE67423} & 150 & \texttt{624C6D4B4351EE0A34728E6B7B6056085064E3} \\
118 & \texttt{60EF5650D1DE83A1E9D48981463723} &  &  \\
\end{longtable}
}

The Newman polynomials obtained are available both in hexadecimal representation and as lists of coefficients in the files \texttt{resulthexa.txt} and \texttt{resultcoeffs.txt}, which can be found at \cite{ComputationalResults}. For example, two entries of Table~\ref{tablehexa}, namely
\[
\texttt{18995E6FE1AA04687E78537E4B41F9B}
\]
and
\[
\texttt{624C6D4B4351EE0A34728E6B7B6056085064E3},
\]
correspond to the following Newman polynomials:

\smallskip

$x^{120}+x^{119}+x^{115}+x^{112}+x^{111}+x^{108}+x^{106}+x^{104}+x^{103}
{}+x^{102}+x^{101}+x^{98}+x^{97}+x^{95}+x^{94}+x^{93}+x^{92}+x^{91}+x^{90}
{}+x^{89}+x^{84}+x^{83}+x^{81}+x^{79}+x^{77}+x^{70}+x^{66}+x^{65}+x^{63}
{}+x^{58}+x^{57}+x^{56}+x^{55}+x^{54}+x^{53}+x^{50}+x^{49}+x^{48}+x^{47}
{}+x^{42}+x^{40}+x^{37}+x^{36}+x^{34}+x^{33}+x^{32}+x^{31}+x^{30}+x^{29}
{}+x^{26}+x^{23}+x^{21}+x^{20}+x^{18}+x^{12}+x^{11}+x^{10}+x^{9}+x^{8}
{}+x^{7}+x^{4}+x^{3}+x+1$,

\smallskip

and 

\smallskip

$x^{150}+x^{149}+x^{145}+x^{142}+x^{139}+x^{138}+x^{134}+x^{133}+x^{131}
{}+x^{130}+x^{128}+x^{126}+x^{123}+x^{121}+x^{120}+x^{118}+x^{113}+x^{112}
{}+x^{110}+x^{108}+x^{104}+x^{103}+x^{102}+x^{101}+x^{99}+x^{98}+x^{97}
{}+x^{91}+x^{89}+x^{85}+x^{84}+x^{82}+x^{78}+x^{77}+x^{76}+x^{73}+x^{71}
{}+x^{67}+x^{66}+x^{65}+x^{62}+x^{61}+x^{59}+x^{57}+x^{56}+x^{54}+x^{53}
{}+x^{52}+x^{51}+x^{49}+x^{48}+x^{46}+x^{45}+x^{38}+x^{36}+x^{34}+x^{33}
{}+x^{27}+x^{22}+x^{20}+x^{14}+x^{13}+x^{10}+x^{7}+x^{6}+x^{5}+x+1$.

\smallskip

For an independent cross-check, we also verified using GP/PARI that these polynomials are indeed Newman multiples of $l(x)^2$.

A natural next step after obtaining these multiples of $l(x)^2$ is to apply the same method in order to search for multiples of $l(x)^3$ within the set $\mathcal{N}$ of Newman polynomials. In \cite{Mossinghoff2003}, M.~J.~Mossinghoff verified that $l(x)^3$ does not divide any Newman polynomial of degree less than $120$. Using our method, we verify that $l(x)^3$ does not divide any Newman polynomial of degree at most $160$.

\section{Conclusion}
The results obtained in this work naturally lead to several open directions for further investigation. In particular, the existence of Newman multiples for the three unresolved polynomials of Mahler measure less than $1.3$ is still unknown, and it remains unclear whether a universal constant governing the existence of Newman multiples for integer-coefficient polynomials according to their Mahler measure exists. From a computational perspective, it is likely that some of the bounds and constructions obtained in this work could be further improved by testing alternative optimization solvers or software environments beyond the one employed here. More generally, the mixed-integer linear programming formulations developed in this article appear very promising and should provide a flexible framework for addressing a wider range of problems involving polynomials with restricted coefficients, including questions related to Borwein, Littlewood and Newman polynomials.

\end{document}